\documentclass{IEEEtran}
\usepackage{cite}
\usepackage{amsmath,amssymb,amsfonts}
\usepackage{algorithmic}
\usepackage{graphicx}
\usepackage{color}
\usepackage{multirow}
\usepackage{textcomp}
\def\BibTeX{{\rm B\kern-.05em{\sc i\kern-.025em b}\kern-.08em
    T\kern-.1667em\lower.7ex\hbox{E}\kern-.125emX}}
\begin{document}
\title{Stability-constrained Optimization for Nonlinear Systems based on Convex Lyapunov Functions}
\author{Qifeng~Li,~\IEEEmembership{Member,~IEEE,}
         and Konstantin~Turitsyn,~\IEEEmembership{Member,~IEEE}
\thanks{The authors are with the Department
of Mechanical Engineering, Massachusetts Institute of Technology, Cambridge
MA, 02139 USA; e-mail: \{qifengli,turitsyn\}@mit.edu }
\thanks{Q. Li is currently with the Department of Electrical $\&$ Computer Engineering, University of Central Florida, Orlando FL, e-mail: qifeng.li@ucf.edu.}}

\maketitle

\begin{abstract}
This paper presents a novel scalable framework to solve the optimization of a nonlinear system with differential algebraic equation (DAE) constraints that enforce the asymptotic stability of the underlying dynamic model with respect to certain disturbances. Existing solution approaches to analogous DAE-constrained problems are based on discretization of DAE system into a large set of nonlinear algebraic equations representing the time-marching schemes. These approaches are not scalable to large size models. The proposed framework, based on \textit{LaSalle'}\textit{s invariance principle}, uses convex Lyapunov functions to develop a novel stability certificate which consists of a limited number of algebraic constraints. We develop specific algorithms for two major types of nonlinearities, namely Lur'e, and quasi-polynomial systems. Quadratic and convex-sum-of-square Lyapunov functions are constructed for the Lur'e-type and quasi-polynomial systems respectively. A numerical experiment is performed on a 3-generator power network to obtain a solution for transient-stability-constrained optimal power flow.

\end{abstract}

\begin{IEEEkeywords}
Convex Lyapunov function, Lur'e systems, quasi-polynomial systems, SOS-convex, Stability-constrained optimization
\end{IEEEkeywords}

\section{Introduction}
\label{sec:introduction}

\subsection{Background and Motivation}
Efficiency and security are often inconsistent for many nonlinear control systems \cite{ZengOptimal}. For instance, improving the stability of power systems may result in a higher operational cost. On the other hand, pursuing the cost-efficiency of power generators may deteriorate the system's security \cite{Gan}. In system identification, stability of the identified model need to be guaranteed while pursuing the minimum modeling errors \cite{LacySubspace}. Efficiency here may represent cost-reduction or improvement in system's performance with respect to some criteria, such as reducing modeling errors in system identification. Similar issue is encountered in many other control systems, such as networked DC motor systems \cite{ZengOptimal}, cognitive radio networks \cite{AlthunibatOn}, and cloud computing/storage \cite{LiSecurity}.  Therefore, it is desirable for the designers or operators of a nonlinear system to achieve a trade-off between efficiency and security.

The dynamics of a nonlinear systems is generally formulated as a set of differential algebraic equations (DAEs). Thus, the security and stability of the system can be quantitatively expressed as the boundedness or convergence of its dynamic trajectories\cite{Khalil}. A straight-forward idea of achieving the optimal trade-off of efficiency and security is to consider the DAEs or partial differential equations (PDEs) as constraints in the optimization problem of maximizing the efficiency. This implementation produces a DAE- or PDE-constrained optimization problem which is not directly solvable by mature optimization algorithms. The most commonly used solution methods for such optimization problems are the discretize-then-optimize approaches \cite{BettsDiscretize,BettsDirect,BettsPractical}. Namely, the basic idea of these methods is to discretize a DAE into a set of algebraic equations with respect to small time steps. 

However, even a small DAE-constrained problem induces a large-scale nonlinear programming (NLP) problem after discretization \cite{BieglerLarge,Scala,Mak,Hijazi}. Moreover, there may be convergence issues associated with the discretize-then-optimize methods \cite{HagerRunge}. As a result, the discretization scheme for DAEs or PDEs needs to be carefully chosen. Therefore, these methods are not practical for large-scale systems like modern transmission grids. A more scalable solution method is desirable for the stability-constrained optimization problems of large-scale nonlinear systems.


\subsection{Relevant Work and Novelty}

This paper develops a novel solution framework for stability-constrained optimization problems in large-scale nonlinear systems and introduces its applications in two typical types of nonlinear systems--Lur'e and quasi-polynomial systems.

Different from the existing methods that discretize the DAEs into a large number of nonlinear algebraic equations, the proposed approach constructs the stability certificate of a nonlinear system based on \textit{LaSalle's Invariance} principle\cite{Khalil}. By utilizing the computationally effective convex form of Lyapunov functions\cite{BoydConvex}, we convert the stability certificate from a set of DAEs into a limited number of algebraic constraints which can be directly incoporated into the optimization framework. As a result, the DEA-constrained optimization is converted into a scalable NLP problem which is more computationally tractable.

We first apply the proposed approach to develop a stability-constrained optimization framework for Lur'e-type systems\cite{WadaParametric}. Quadratic Lyapunov function candidates, which are convex, can be constructed for a Lur'e system with a properly chosen sector \cite{Sector} for the nonlinear terms\cite{BoydLinear}. A swing equation based dynamic model of a power transmission network was reformulated into the Lur'e form in \cite{LongA}. In the numerical experiment section, we demonstrate the proposed approach in a power transmission network \cite{Kundur} whose dynamic model is a Lur'e-type system. 

The applicability of the proposed approach to a type of quasi-polynomial systems\cite{PapachristodoulouAnalysis} is also investigated in this paper. The quasi-polynomial systems represent a much wider range of nonlinear systems than the Lur'e-type systems do. The quasi-polynomial system is first re-cast into a polynomial system. Then, a convex sum of squares (SOS-convex) Lyapunov function \cite{AhmadiSOS,KunduA} is constructed for this system. The dynamic model of a microgrid may not be of Lur'e-form\cite{VorobevHigh}. However, almost all existing microgrid dynamic models belong to the class of quasi-polynomial systems, implying that the proposed approach is also valid for achieving the optimal trade-off of security and efficiency in microgrids.

The rest of the paper is organized as follows: Section II introduces the mathematical formulation of the problem. The novel stability-constrained optimization framework for nonlinear systems is presented in Section III. The applications in Lur'e and quasi-polynomial systems are introduce in Section IV, where the methods for constructing convex Lyapunov functions are also investigated. In Section V, a numerical experiment is given for validating the proposed approach, while conclusion is drawn in Section VI.

\section{Problem Formulation and Statement}
\subsection{A Dynamically-constrained Optimization Problem}
\begin{figure}
\centering
\includegraphics[width=0.48\textwidth]{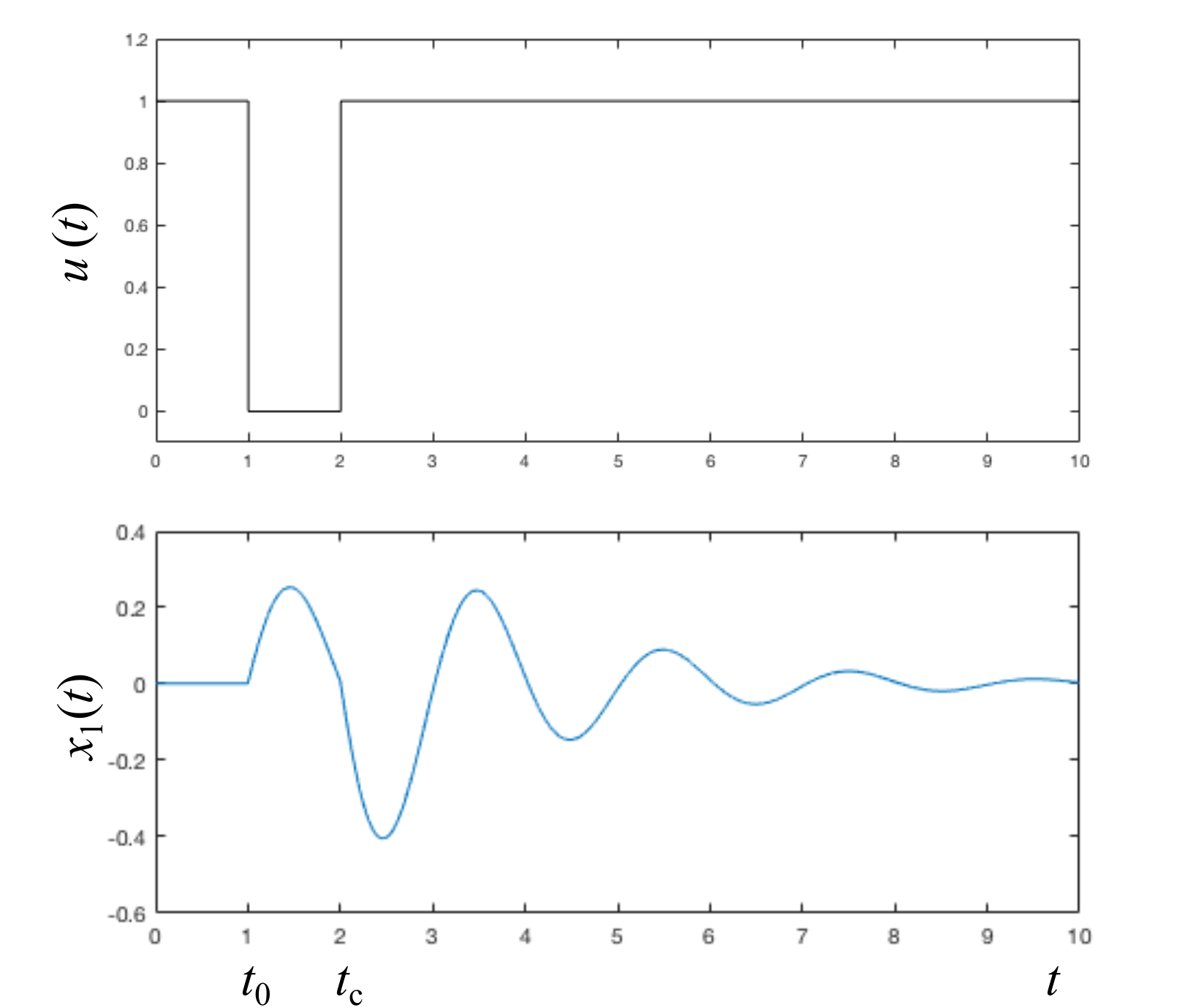}
\caption{An example of disturbance $u(t)$ and its impact on the system trajectory of a modified pendulum system whose formulation is give in Appendix \ref{app:pendulum}.}
\label{fig:pulse}
\end{figure}

This paper considers the following DAE-constrained optimization problem:
\begin{subequations} \label{o-SCO}
\begin{align}
\min_{w} \;\, &c(w) \label{Objective}  \\ 
\mathrm{s.t.}\;\, & \dot{x}=f(x,w,u(t)) \label{DynamicSyst}\\
&0=f(x^o,w,0) \label{NonLinearSyst} \\
&0 \ge g(x^o,w) \label{PhysicalConstraint} \\
& \|x(t) - x^o \| \le \epsilon \; \text{as}\; t  \rightarrow +\infty, \label{Converge}
\end{align}
\end{subequations}
where $x^o \in \mathbb{R}^n$ and $x\in \mathbb{R}^n$ are column vectors of steady-state and dynamic-state variables respectively, $\epsilon$ is an arbitrarily small positive scalar, and $x^o=x|_{t =t_0}$. $c(w):\mathbb{R}^m \rightarrow \mathbb{R}$ is a cost or performance function of the control variable vector $w$. An example of the vanishing disturbance $u(t) \in \mathbb{R}^r$ and its impact on the system trajectory $x(t)$ is given in Figure \ref{fig:pulse}. $f:\mathbb{R}^{n+m+r} \rightarrow \mathbb{R}^n$ is a column vector of nonlinear functions. Generally, inequality (\ref{PhysicalConstraint}) represents the bound constraints, such as $\underline{x}^o \le x^o \le \overline{x}^o$ and $\underline{w} \le w \le \overline{w}$. An engineering example of problem (\ref{o-SCO}) is provided in Section V.

Suppose that the nominal system (i.e. the dynamic system without disturbance) of (\ref{DynamicSyst}):
\begin{equation} \label{nominal}
\dot{x}=f(x,w,0),
\end{equation}
is locally stable in a neighbourhood $\mathcal{P}$ of the equilibrium $x^o$, the stability of the perturbed system is jointly determined by $u$ (including the norm and duration of $u$) and $x^o$. The goal of solving problem (\ref{o-SCO}) is to select the optimal solution of $w$ with regard to its cost function $c(w)$ and, at the same time, guarantee that the resulting steady-state operating point $x^o$ can withstand a given disturbance $u$.

The DAE-constrained optimization (\ref{o-SCO}) is not directly solvable by the existing optimization algorithms. An existing approach differentiates the differential equation (\ref{DynamicSyst}) into multiple algebraic equations resulting in a nonlinear program (NLP) with a large number of nonlinear constraints. This approach is not scalable for large-scale nonlinear control system, such as power grids. This paper aims at developing a computationally tractable solution method to problem (\ref{o-SCO}) based on the Lyapunov stability theories.

\section{Stability-constrained Optimization}
This section presents a novel scalable framework for solving the dynamically-constrained optimization problem (\ref{o-SCO}). The primary idea of the proposed approach is to convert the stability condition (\ref{DynamicSyst}) and (\ref{Converge}) into a limited number of algebraic constraints based on \textit{LaSalle's Invariance} principle, so that the resulting optimization model is an ordinary NLP problem which can be directly solved by mature algorithms or solvers, such as KNITRO\cite{knitro} and IPOPT \cite{ipopt}.

\subsection{Scalable Algebraic Stability Certificate}
Let
\begin{equation}
\mathcal{F}=\left\{(x^o,w) \in \mathbb{R}^{n+m} \left| \begin{array}{lr}
f(x^o,w,0)=0\\
g(x^o,w) \le 0 
\end{array}\right. \right\},\nonumber
\end{equation}
we have the following definition:

\textbf{Definition 1.} \textit{A function} $V(x):\mathbb{R}^n \rightarrow \mathbb{R}$ \textit{is said to be a candidate of common local Lyapunov functions (CLLFs) for the dynamic system (\ref{nominal}) if the following conditions hold:} \\
\textit{a)} there exists a polytope $\mathcal{P}$ for each $x^o \in \mathcal{F}$ that
\begin{equation}
x^o \in \mathcal{P} =\{x \in \mathbb{R}^n \mid Cx -d\le 0\} \nonumber
\end{equation}
\textit{where $C$ is a $|\mathcal{N}| \times n$ matrix while $d$ is a $|\mathcal{N}|$-dimensional vector, $\mathcal{N}$ is the facet set of $\mathcal{P}$;}\\
\textit{b)} $V(x^o)=0$ \textit{and} $V(x)>0$ ($\forall x \in \mathcal{P}-\{x^o\}$ \textit{and} $\forall (x^o,w) \in \mathcal{F}$)\textit{;}\\
\textit{c) } $\mathcal{P} \subset \mathcal{D}=\{x \in \mathbb{R}^n \mid \dot{V}(x) \leq 0\}$ \textit{for all} $(x^o,w) \in \mathcal{F}$.


Suppose there exists a \textit{CLLF}, $V(x)$, for the nominal system (\ref{nominal}), we have the following lemma of \textit{LaSalle'}\textit{s Invariance} principle:

\textbf{Lemma 1}. \textit{The system trajectories of (\ref{nominal}) starting from a fault-cleared state $x^c=x|_{t=t_c}$ stay in the set $\mathcal{R}$ for all $t \in [t_c,+\infty)$ and eventually will converge to $x^o$ if:}\\
\textit{a)} the fault-cleared state $x^c$ is
within the set $\mathcal{R}$ which is defined as
\begin{equation} \label{Vmin}
 \mathcal{R} = \{x \in \mathcal{P} \mid V(x) \le V^{\min}\},
  \end{equation}
\textit{and}
\begin{align} \label{Vmindefinition}
V^{\min}=\min_{x \in \partial \mathcal{P}}\;&V(x) ,
\end{align}
\textit{where the boundary of} $\mathcal{P}$
\begin{equation}
\partial \mathcal{P}= \bigcup_{i \in \mathcal{N}} \partial \mathcal{P}_i =\bigcup_{i \in \mathcal{N}} \left\{\textbf{\emph{x}} \left| \begin{array}{lr}
C_i^{\mathrm{T}}x -d_i= 0 \\
Cx-d \le 0 ,  
\end{array}\right. \right\}, \nonumber
\end{equation}
\textit{$C_i^{\mathrm{T}}$ is the $i$th row of $C$, and $d_i$ is the $i$th element of vector $d$;}\\
\textit{b)} $V(x)$ \textit{is a valid CLLF for system (\ref{nominal}).}


\textit{Proof}: See Appendix \ref{app:lemma1}. $\square$
\begin{figure}
\centering
\includegraphics[width=0.3\textwidth]{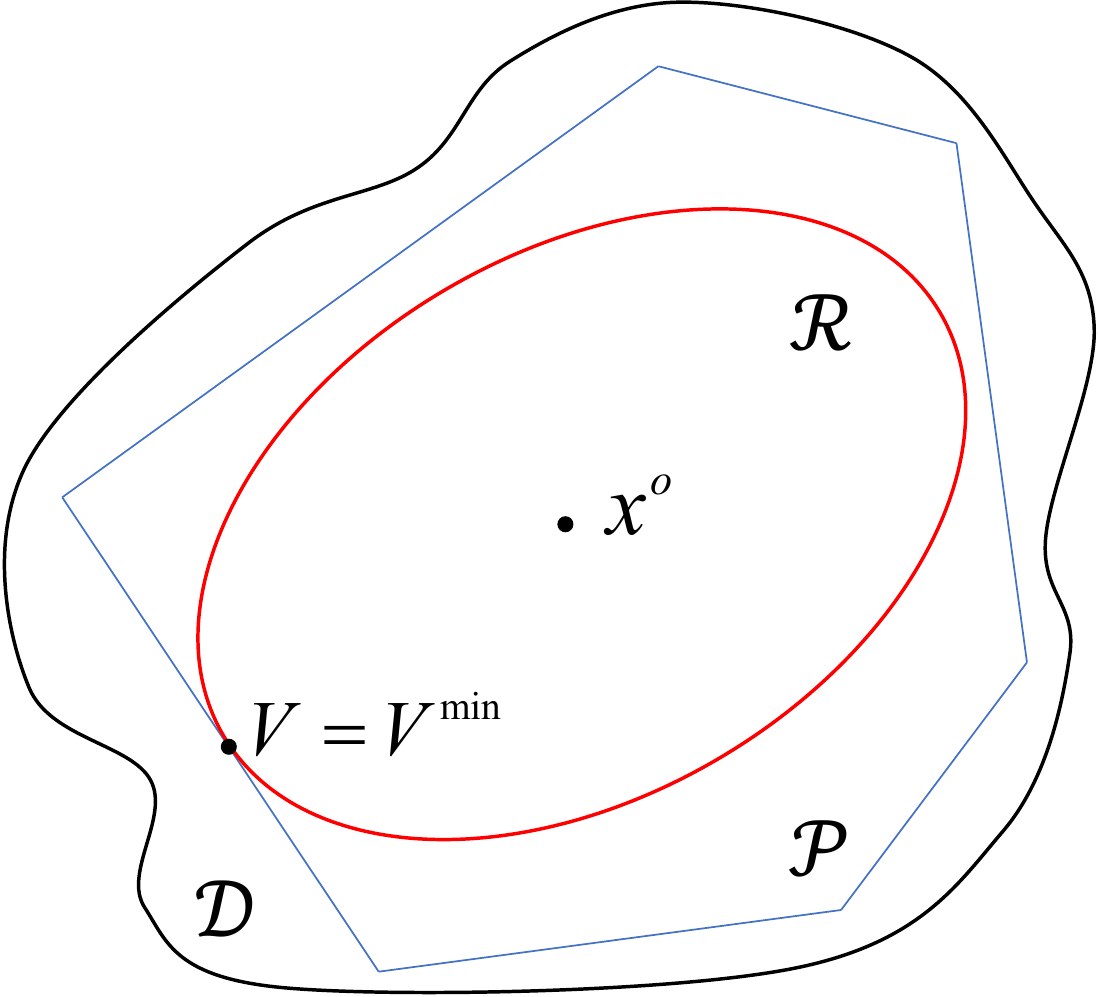}
\caption{Geometric representation of sets in Definition 1 and Lemma 1.}
\label{fig:sets}
\end{figure}

A geometric representation of sets $\mathcal{P}$, $\mathcal{D}$, and $\mathcal{R}$ is given in Figure \ref{fig:sets}. The basic idea of Definition 1 and Lemma 1 is to search for the minimum Lyapunov value along the facets of a pre-selected polytope $\mathcal{P}$ inside $\mathcal{D}$ rather than on the boundary of $\mathcal{D}$.  The notation $\partial \mathcal{P}$ represents the boundary of $\mathcal{P}$ that can be formulated as a disjunctive set which is the union of disjunctions $\partial \mathcal{P}_i$ ($i=1,\,\ldots,\, |\mathcal{N}|$). In reality, $\partial \mathcal{P}_i$ represents a facet of the polytope $\mathcal{P}$ from the perspective of geometry. Consequently, optimization (\ref{Vmindefinition}) is a disjunctive programming (DP) which can not be directly solved by mature optimization algorithms. To mitigate this issue, we propose the following NLP reformulation of (\ref{Vmindefinition}):
\begin{subequations} \label{newVmin}
\begin{align} 
V^{\min}=\min_{x^i}\;& V(x^i) \\
\mathrm{s.t.}\; &\frac{\partial V(x^i)}{\partial x^i} - \lambda_i C_i =0 \\
&C_i^{\mathrm{T}}x^i -d_i= 0 \\
&Cx^i-d \le 0 , 
\end{align}
\end{subequations}
where $x^i \in \mathbb{R}^n$ ($i \in \mathcal{N}$) is the state vector on the $i$th facet of the polytope $\mathcal{P}$ and $\lambda_i \in \mathbb{R}$ is the $i$th Lagrange multiplier ($i \in \mathcal{N}$).  Then, we have the following theorem.

\textbf{Theorem 1}. \textit{The $V^{\min}$ obtained by solving optimization problem}  (\ref{newVmin}) \textit{is the same as the optimal solution of} (\ref{Vmindefinition}) \textit{if $V(\cdot)$ is convex.}

\textit{Proof}: See Appendix \ref{app:theorem1}. $\square$

\subsection{A Scalable Stability-constrained Optimization Framework}
At $t=t_0$, system (\ref{DynamicSyst}) suffers from a disturbance and its trajectory $x(t)$ starts deviating from the equilibrium point $x^o$. In reality, the fault-cleared state $x^c$ is a function of $x^o$ if the disturbance $u$ and the fault clearing time $t_c$ is given. Via Taylor's series, such a function can be explicitly expressed as
\begin{align} \label{trajectory}
  x^c = x^o +\sum_{n=1}^{N} \frac{f^{(n)}(x^o,w,u^o)}{n!}(t_c-t_0)^n,
\end{align}
where $f^{(n)}$ is the $(n-1)$th-ordered derivative of $f$ with respect to time $t$; $u^o=u|_{t=t_0}$ and $N=\infty$. A small $N$ can provide satisfactory accuracy if ($t_c-t_0$) is small.

According to Lemma 1, the post-disturbance trajectory $x(t \in [t_c,+\infty))$ can stay within $\mathcal{P}$ or even converge back to $x^o$ as $t \rightarrow + \infty$ if the fault-clearing state $x^c$ is in $\mathcal{R}$, namely
\begin{equation}
V(x^c) \le V^{\min}. \label{energy}
\end{equation}
Note that $x^c$ in (\ref{trajectory}) is an explicit function, while $V^{\min}$ in (\ref{newVmin}) is a nonexplicit function, of $x^o$. Hence, constraints (\ref{newVmin})-(\ref{energy}) together represent the dynamic stability certificate in the $x^o$-domain. By replacing constraints (\ref{DynamicSyst}) and (\ref{Converge}) with (\ref{newVmin})-(\ref{energy}), we have the following stability-constrained optimization model for the nonlinear system (\ref{DynamicSyst}):
\begin{align}  
\begin{split}
\min_{w} \quad &\text{(\ref{Objective})} \\
\text{s.t.} \quad &\text{(\ref{NonLinearSyst}),\,(\ref{PhysicalConstraint})},\,\text{and}\,(\ref{newVmin})-(\ref{energy}).
\end{split} \tag{b-SCO}
\end{align}

We can consider (b-SCO) as an alternative formulation to the DAE-constrained optimization model (\ref{o-SCO}). However, model (b-SCO) is still very hard to solve since it is a bilevel optimization problem where (\ref{newVmin}) is the lower-level subproblem. To overcome this issue, we proposed the following single-level NLP model
\begin{subequations}
\begin{align}  
\text{(s-SCO)}\quad\min_{w,V^{\min}}\; &c^\prime(w,V^{\min}) = c(w) - \epsilon V^{\min} \label{Objective1}\\
\text{s.t.} \; &\text{(\ref{NonLinearSyst}),\,(\ref{PhysicalConstraint})},\,(\ref{trajectory}),\,(\ref{energy}),\,\text{(5b)-(5d)},\,\text{and} \nonumber \\
&V^{\min} \le V(x^i), \label{Vminbound} 
\end{align}
\end{subequations}
where $\epsilon$ is an arbitrarily small positive value and $i \in \mathcal{N}$.


\textbf{Theorem 2}. \textit{Optimal solution of} (s-SCO) \textit{is also optimal to optimization problem} (b-SCO).

\textit{Proof}: See Appendix \ref{app:theorem2}. $\square$

\textbf{Remark 1}. Theorem 1 demonstrates that the single-level NLP problem (s-SCO) is equivalent to the bilevel problem (b-SCO), which implies (s-SCO) is also a stability-constrained optimization framework for nonlinear system (\ref{DynamicSyst}). Solving the proposed SCO framework (s-SCO), which is a scalable NLP problem, to obtain an optimal solution to  (\ref{o-SCO}) outperforms the discrete-then-optimize approach in terms of computational efficiency due to less number of nonlinear constraints.



\section{Construction of Convex Lyapunov Functions}
The analysis in the previous section is based on an assumption that there exists a convex \textit{CLLF} for the nonlinear system. This section introduces the application of the proposed approach in two types of nonlinear systems--Lur'e and  quasi-polynomial systems--together with the methods of constructing convex \textit{CLLF}s.

\subsection{Lur'e-type Systems}
If the dynamic system (\ref{nominal}) can be rewritten as
\begin{equation} \label{Lure}
\dot{x}=f(x,w,0)=A(x-x^o)+B\phi(C(x-x^o)),
\end{equation}
where $\phi_i:\mathbb{R} \rightarrow \mathbb{R}$ $(i=1,\ldots,l)$ is nonlinear function, and $A$ and $B$ are $n \times n$ and $n \times l$ matrices respectively, it is a Lur'e-type system. Equations (\ref{DynamicSyst}) and (\ref{NonLinearSyst}) imply that $x^o$, $A$, and $\phi$ are related to the vector of parameters $w$.

\textbf{Definition 2.} \textit{The nonlinear function $\phi(x)$ is said to be locally bounded by sector} [$\gamma,\beta$] \textit{if, for all $x \in \mathcal{P}$ and $(x^o,w) \in \mathcal{F}$,} 
\begin{equation} \label{Sector}
(\phi- \gamma C(x-x^o))^\mathrm{T}(\phi - \beta C(x-x^o)) < 0 .
\end{equation}

\begin{figure}
\centering
\includegraphics[width=0.32\textwidth]{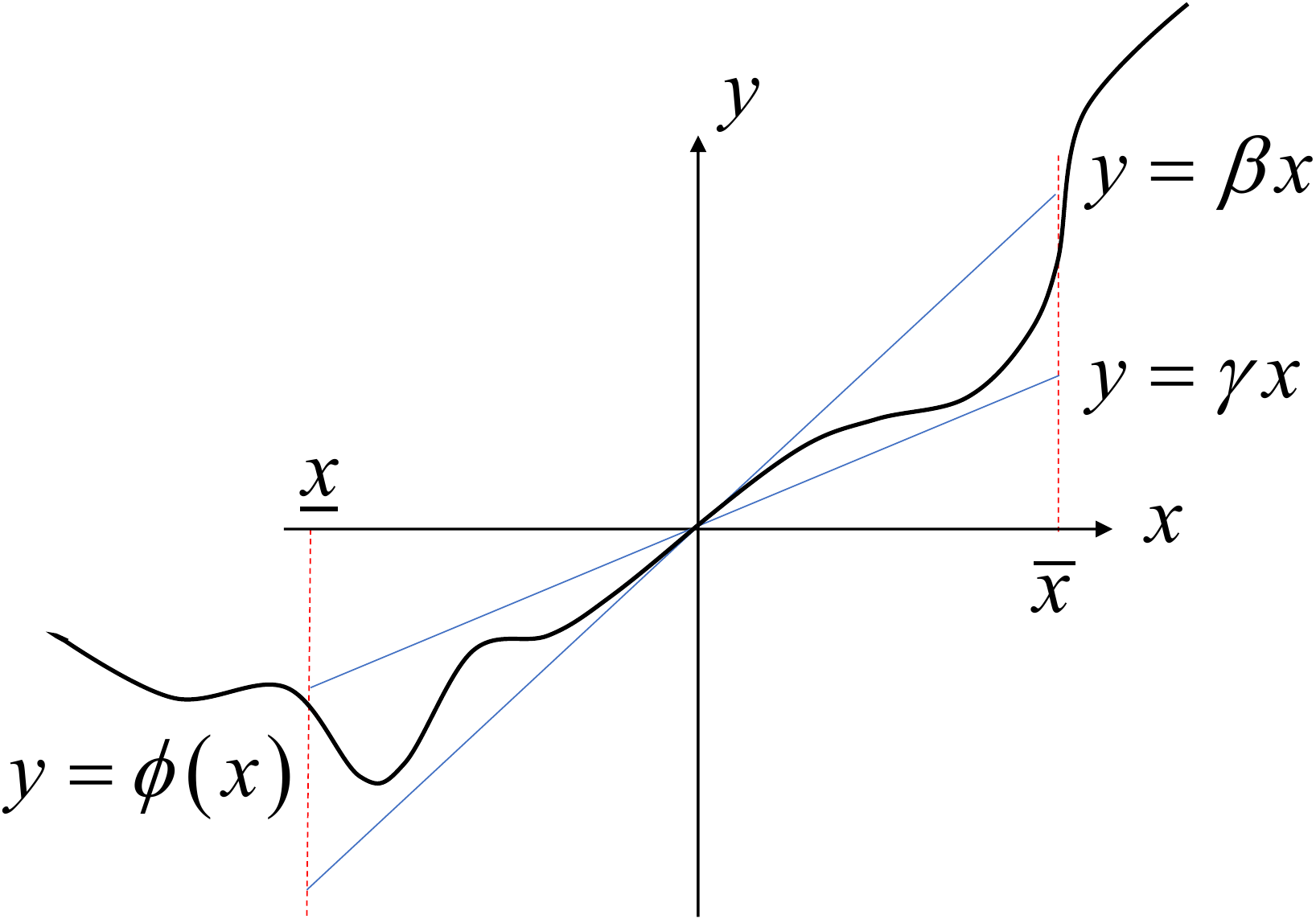}
\caption{Geometric interpretation of sector bounds.}
\label{fig:sectorbounds}
\end{figure}

A geometric interpretation of sector bounds is given in Figure \ref{fig:sectorbounds}, where the upper and lower bounds of $x$ is corresponding to the facets of $\mathcal{P}$. Generally, the sector [$\gamma,\beta$] and facets of $\mathcal{P}$ are synergistically determined. For Lur'e systems, we directly use the matrix $C$ in the formulation (\ref{Lure}) to construct the following polytope with parallel facets:
\begin{equation}
\mathcal{P} =\{x \in \mathbb{R}^n \mid \underline{d} \le C^\mathrm{T}x \le \overline{d}\}. \nonumber
\end{equation}

Based on Definition 2, we can construct the following convex quadratic \textit{CLLF}
\begin{equation} \label{QuaLyapunov}
 V(x)=(x-x^o)^\mathrm{T}P(x-x^o)
\end{equation}
for system (\ref{Lure}) by solving the following linear matrix inequality (LMI) \cite{BoydLinear} for the positive matrix $P$
\begin{equation} \label{LMI}
\left[
\begin{array}{cc} A^\mathrm{T}P+PA -   C^\mathrm{T}\tau\gamma\beta C     & PB+\frac{1}{2}(\gamma+\beta)C^\mathrm{T} \\
 B^\mathrm{T}P+\frac{1}{2}(\gamma+\beta)C  & -\tau \\
\end{array}
\right] \preceq 0.
\end{equation}
where $\tau \geq 0$, if $A$ is Hurwitz. As discussed in the proof of Theorem 1, the optimal solution of $x^i$ is the unique solution $\hat{x}^i$ determined by (5b)-(5d). Since, for the Lur'e cases, the constructed \textit{CLLF} $V$ is a quadratic function, $\hat{x}^i$ ($i \in \mathcal{N}$) has a non-trivial expression 
\begin{equation} 
\hat{x}^i=\frac{\underline{d}_iP^{-1} C_i}{C_i^\mathrm{T}
P^{-1} C_i} \;\text{or}\;\frac{\overline{d}_iP^{-1} C_i}{C_i^\mathrm{T}
P^{-1} C_i}. \nonumber
\end{equation}
As a result, the corresponding Lyapunov value $V(\hat{x}^i)$ is 
\begin{equation} \label{Concavee}
V_i^{\min}=\frac{\min\left\{(C_i^\mathrm{T}x^o - \underline{d}_i)^2,\;(C_i^\mathrm{T}x^o - \overline{d}_i)^2\right\}}{C_i^\mathrm{T}
P^{-1} C_i} .
\end{equation}
With applying the stability-constrained optimization framework (s-SCO) to the Lur'e cases, we have
\begin{align}  
\text{(SCO1)}\quad\quad\min_{x^o,w}\quad &\text{(\ref{Objective1})} \nonumber \\
\text{s.t.} \quad &\text{(\ref{NonLinearSyst}),\,(\ref{PhysicalConstraint})},\,(\ref{trajectory}),\,(\ref{energy}),\,\text{and} \nonumber \\
\begin{split} \label{Concave}
&V^{\min} \le \frac{(C_i^\mathrm{T}x^o - \underline{d}_i)^2}{C_i^\mathrm{T}
P^{-1} C_i} \\
&V^{\min} \le \frac{(C_i^\mathrm{T}x^o - \overline{d}_i)^2}{C_i^\mathrm{T}
P^{-1} C_i},
\end{split}
\end{align}
where $i \in \mathcal{N}$.

Convex optimization, due to its high-efficiency, has been applied to a wide range of automatic control systems \cite{BoydConvex,LiConvex1}. To facilitate the future development of a convex SCO model for Lur'e system (\ref{Lure}), we discuss two convex versions of the concave constraint (\ref{Concave}) in this section: convex relaxation (\ref{CR}) and convex inner approximation (\ref{CI}). 
\begin{subequations}\label{CR}
\begin{gather}
V^{\min} \le \frac{a_iC_i^\mathrm{T}x^o+b_i}{C_i^\mathrm{T}
P^{-1} C_i}  \\
V^{\min} \le \frac{a_i^\prime C_i^\mathrm{T}x^o+b_i^\prime}{C_i^\mathrm{T}
P^{-1} C_i}. 
\end{gather}
\end{subequations}
\begin{subequations}\label{CI}
\begin{gather}
V^{\min} \le \frac{c_iC_i^\mathrm{T}x^o-\overline{C_i^\mathrm{T}x^o}\underline{C_i^\mathrm{T}x^o}+\Delta d_i^2}{C_i^\mathrm{T}
P^{-1} C_i} \\
V^{\min} \le \frac{2\Delta d_i C_i^\mathrm{T}x^o-\Delta d_i^2}{C_i^\mathrm{T}
P^{-1} C_i}. 
\end{gather}
\end{subequations}
where $a_i=\overline{C_i^\mathrm{T}x^o} -3\overline{d}/2+\underline{d}/2$, $b_i=\overline{d}^2-(\overline{d}+\underline{d})\overline{C_i^\mathrm{T}x^o}/2$, $a_i^\prime=\underline{C_i^\mathrm{T}x^o} -3\underline{d}/2+\overline{d}/2$, $b_i^\prime=\underline{d}^2-(\overline{d}+\underline{d})\underline{C_i^\mathrm{T}x^o}/2$, $c_i=\overline{C_i^\mathrm{T}x^o} +\underline{C_i^\mathrm{T}x^o}+ 2\Delta d_i$and (\ref{CI}) is developed by considering the intersection point of (\ref{Concave}) and the $V^{\min}$-axis as the tangent point. If there exist parallel facets of polytope $\mathcal{P}$, namely $d_{i,1} \le C_i^\mathrm{T}(x-x^o) \le d_{i,2}$, we have the following theorem.


\textbf{Theorem 3}. \textit{Set $\Psi$ is the convex hull of set $\psi$, where}
 \[
    \psi = \left\{(x^o,V^{\min}) \left| \text{(\ref{Concave}) and }\underline{C_i^\mathrm{T}x^o} \le C_i^\mathrm{T}x^o \le \overline{C_i^\mathrm{T}x^o} \right. \right\} 
  \]
   \[   
    \Psi = \left\{(x^o,V^{\min}) \left| \text{(\ref{CR}) and }\underline{C_i^\mathrm{T}x^o} \le C_i^\mathrm{T}x^o \le \overline{C_i^\mathrm{T}x^o}\right. \right\}, 
  \]
and $\underline{C_i^\mathrm{T}x^o} \le (\overline{d}+\underline{d})/2 \le \overline{C_i^\mathrm{T}x^o}$.

\textit{Proof}: See Appendix \ref{app:theorem3}. $\square$
\begin{figure}[tb]
\centering
\includegraphics[width=0.48\textwidth]{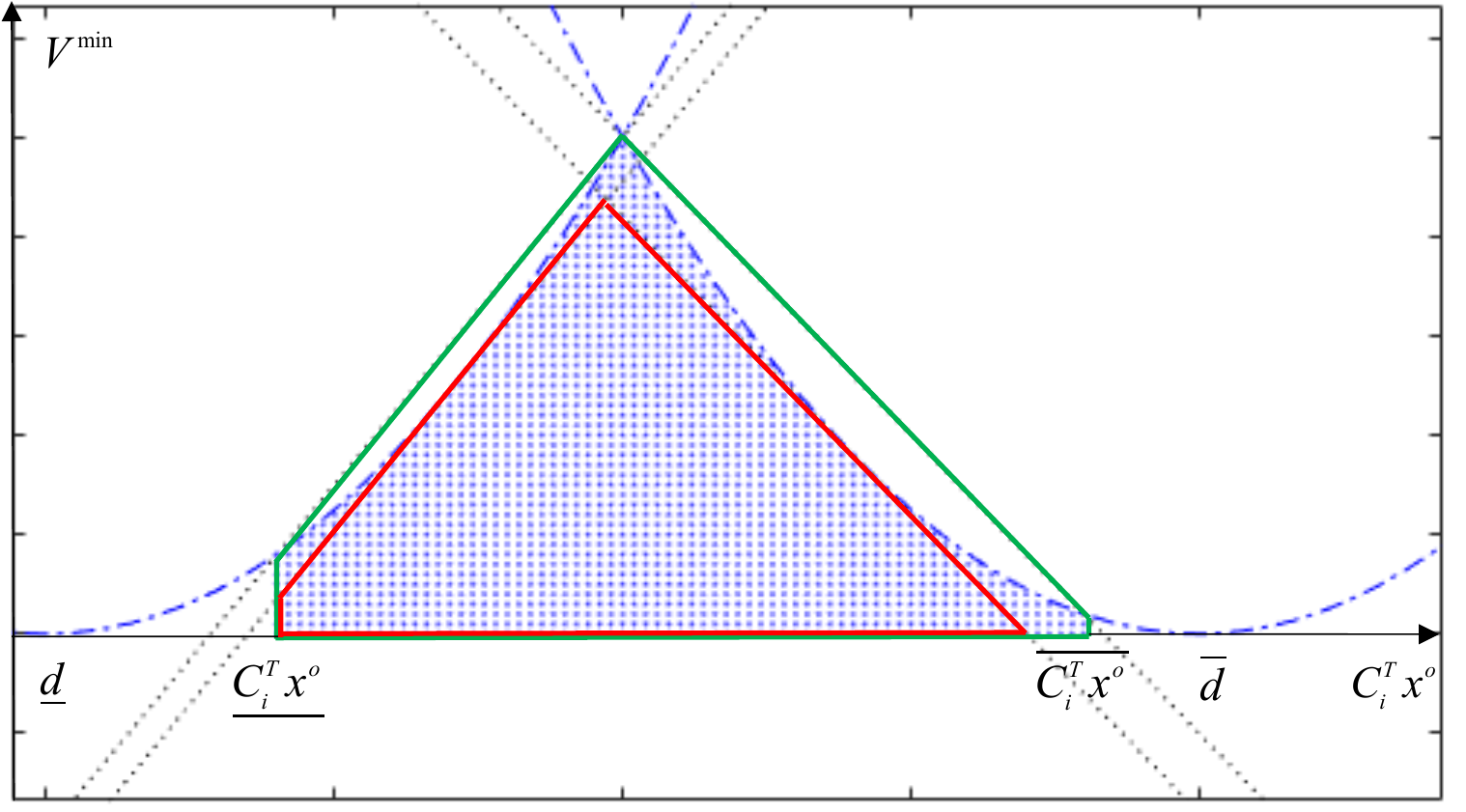}
  \caption{Two convex alternatives of concave constraints (\ref{Concave}). The shaded area denotes the original concave feasible set of (\ref{Concave}). The regions with green and red boundaries are the convex hull relaxation and a convex inner approximation of the shaded area respectively.}
  \label{fig:convexhull}
\end{figure}


It can be observe from figure 2 that, if there exist parallel facets, $\Psi$ is the convex hull\cite{LiConvex2,LiThe}  of the concave constraint (\ref{Concave}) which is tighter than (\ref{CR}). The inner approximation\cite{Hung} can strictly guarantee stability by introducing additional conservativeness. Actually, there exist multiple choices of convex inner approximations of a nonconvex set and the best convex inner approximation is yet to be defined. Defining such an inner approximation will be one of the tasks in our future research. For the cases that quadratic Lyapunov functions are not very conservative, the inner approximation of (\ref{Concave}) is recommended to strictly guarantee stability.  

%



\subsection{Quasi-polynomial Systems}


In this subsection, we consider a nonlinear system in the following form
\begin{equation}   \label{Pseudo}
x_i = f_i(x,w,0)=\sum_j(e_j\prod_k\rho_{ijk}(x-x^o)),
\end{equation}
where $x_i \in \mathbb{R}$ is the $i$th element of $x$, $f_i:\mathbb{R}^{n+m+r} \rightarrow \mathbb{R}$ is the $i$th element of the vector function $f$; $e_j$'s are real numbers, and $\rho_{ijk}$ are elementary functions, or nested elementary functions of elementary functions. Note that $x^o$, $e_j$'s, and $\rho_{ijk}$ are related to $w$ due to constraints (\ref{DynamicSyst}) and (\ref{NonLinearSyst}). System (\ref{Pseudo}) is called  a quasi-polynomial system since it can be converted into a polynomial system whose dimension is higher than that of the original system. For such a purpose, we introduce the following algorithm:

\textbf{Algorithm 1}\footnote{Algorithm 1 is adopted from \cite{PapachristodoulouAnalysis} and \cite{SavageauRecasting} with some modifications. An illustrative example was provided in \cite{PapachristodoulouAnalysis} showing how the algorithm works.}

\begin{enumerate}
\item \textit{For each $\rho_{ijk}(x)$ that is not of the form $\rho_{ijk}(x)=x_l^z$, where $z$ is some integer and $1 \le l \le n$, introduce a new variable $y_m$ and let $y_m=\rho_{ijk}(x)$.}
\item \textit{Compute the differential equation describing the time evolution of $y_m$ using the chain rule of differentiation.}
\item \textit{Replace all appearances of such $\rho_{ijk}(x)$ in the system equations by $y_m$.}
\item \textit{Repeat steps 1)-3), until we obtain system equations of polynomial form.}
\end{enumerate}

With Algorithm 1 applied, system (\ref{Pseudo}) can be recast into the following polynomial system
\begin{subequations} \label{qps}
\begin{gather}
\dot{x}=h_x(x-x^o,y-y^o) \label{qps1}\\
\dot{y}=h_y(x-x^o,y-y^o) \label{qps2} \\
y-y^o=s(x-x^o), \label{qps3}
\end{gather}
\end{subequations}
where constraint (\ref{qps3}) arises directly from the recasting process. Let $\tilde{x}=[x, y]^T$ and $h=[h_x,h_y]^T$, system (\ref{qps}) can be rewritten as
\begin{subequations} \label{Polynomial1}
\begin{gather}
\dot{\tilde{x}}=h(\tilde{x}-\tilde{x}^o)\\
E_2(\tilde{x}-\tilde{x}^o)=s(E_1(\tilde{x}-\tilde{x}^o))
\end{gather}
\end{subequations}
where $E_1$ and $E_2$ are incidence matrices, such that (19b) is equivalent to (\ref{qps3}). We have the following proposition.

\textbf{Proposition 2} \textit{If, for system} (\ref{Polynomial1}) \textit{with all $x^o \in \mathcal{F}$, there exists a polynomial function $V(\tilde{x})$ such that} \\
\textit{a) $V(\tilde{x}^o)=0$;}\\ 
\textit{b) $V(\tilde{x})$ is SOS-convex;}\\
\textit{c) $-\dot{V}(\tilde{x})-p(\tilde{x})$ is SOS, where $p(\tilde{x})$ is a polynomial function which is positive in the set}
\begin{equation}
\mathcal{S}=\{\tilde{x} \mid E_2(\tilde{x}-\tilde{x}^o)-s(E_1(\tilde{x}-\tilde{x}^o))=0\} ;\nonumber
\end{equation}
\textit{d) there exists a polytope} $\mathcal{P}=\{\tilde{x} \mid \tilde{C}\tilde{x} \le \tilde{d}\} \subset \mathcal{S}$\textit{;} \\
\textit{$V(\tilde{x})$ is a valid CLLF for system} (\ref{Polynomial1}).

\textit{Proof}: The SOS-convex Lyapunov function $V$ introduced in the above proposition satisfies all the conditions in Definition 1. Therefore, $V$ is a \textit{CLLF} for the polynomial system (\ref{Polynomial1}). $\square$

To obtain a valid SOS-convex \textit{CLLF} for polynomial system (\ref{Polynomial1}), we introduce the following proposition. For a given polynomial system like (\ref{Polynomial1}), the Matlab toolbox SOSTOOLS\cite{PrajnaIntroducing} can automatically  produce a SOS Lyapunov function, which is not guaranteed to be convex, if there exists one. According to Proposition 3, the SOSTOOLS can search for an SOS-convex \textit{CLLF} if condition (\ref{sos-convex}) is incorporated.

\textbf{Proposition 3}. \textit{Let the polynomial $V(x)$ be SOS, it is also SOS-convex if}
\begin{equation} \label{sos-convex}
(1-\alpha)V(x)+\alpha V(y)-V((1-\alpha)x+\alpha y) \; is \; SOS,
\end{equation}
\textit{where $0 < \alpha < 1$.}

\textit{Proof}: Condition (\ref{sos-convex}) implies that
\begin{equation}
(1-\alpha)V(x)+\alpha V(y)-V((1-\alpha)x+\alpha y) \ge 0, \nonumber
\end{equation}
which means $V(\cdot)$ is a convex function. $\square$

With the SOS-convex \textit{CLLF} obtained above, we can develop the following stability-constrained optimization model for the quasi-polynomial system (\ref{Pseudo}) by applying the framework (s-SCO):

\begin{align}  \tag{SCO2}
\begin{split}
\min_{x^o,w}\quad &\text{(\ref{Objective1})}  \\
\text{s.t.} \quad &\text{(\ref{NonLinearSyst}),\,(\ref{PhysicalConstraint})},\,(\ref{trajectory}),\,\text{and}  \\
&y^c-y^o=s(x^c-x^o) \\
&V(\tilde{x}^c) \le V^{\min}\\
&V^{\min} \le V(\tilde{x}^i) \\
&\frac{\partial V(\tilde{x}^i)}{\partial \tilde{x}^i} - \lambda_i c_i^T =0 \\
&C_i^\mathrm{T}\tilde{x}^i - d_i=0\\
&C\tilde{x}^i - d \le 0, 
\end{split}
\end{align} 
where $\tilde{x}^c=[x^c, y^c]^T$.

In this paper, the set of polynomial systems is considered as a subset of quasi-polynomials. For the cases of polynomial, there is an easier way to implement the SCO framework presented in Section III, which is introduced in Appendix \ref{app:polynomial}.

\section{Application and Experimental Results}

\begin{figure}[h]
\centering
\includegraphics[width=0.45\textwidth]{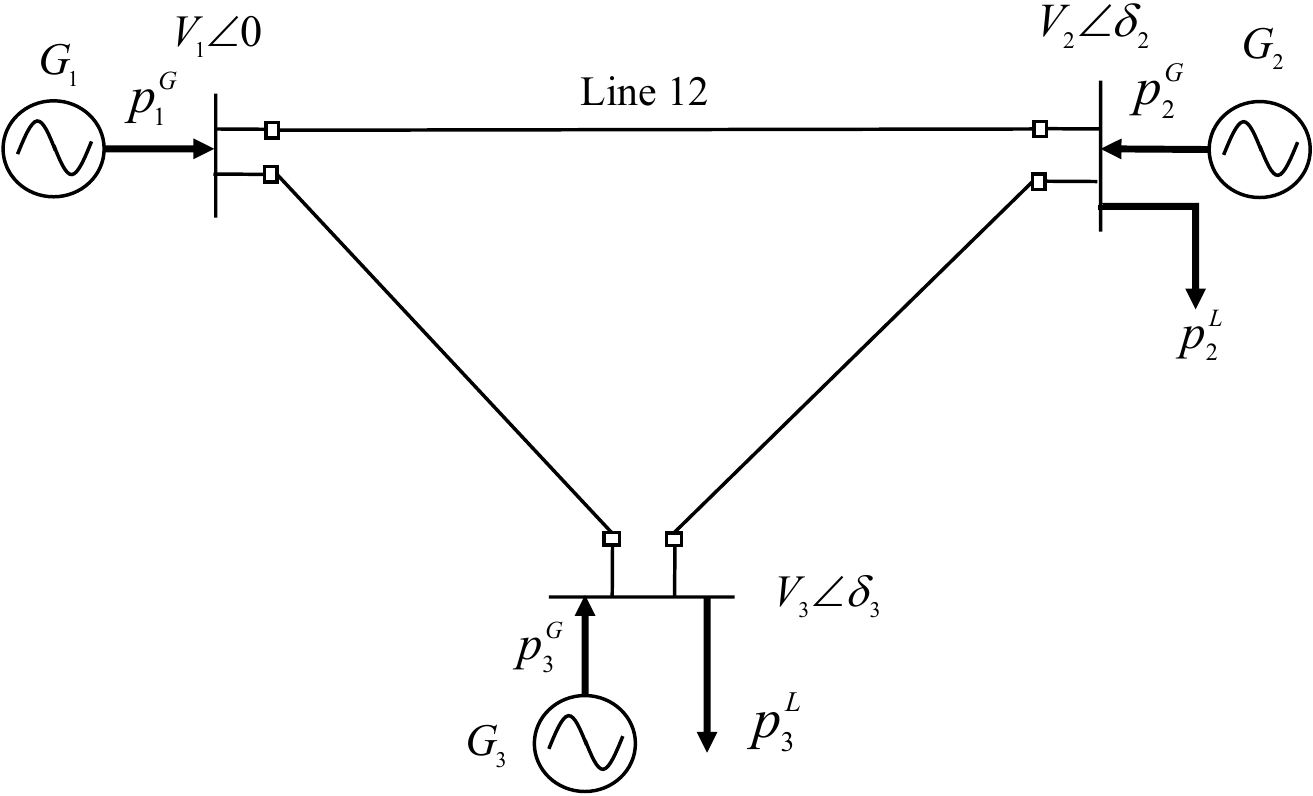}
\caption{The three-generator system.}
\label{fig:3bus}
\end{figure}

This section demonstrates the proposed stability-constrained optimization framework on a three-bus power system shown in Figure \ref{fig:3bus} to obtain a low-cost generation scheme for power grids with transient stability guaranteed. 

\subsection{Problem Statement}
The operators of power grids need to dispatch the generators to meet the load demands through the transmission network with a possibly low cost. Such a low-cost generation scheme can be achieved by solving the following optimal power flow (OFF) problem\cite{OPF}:
\begin{subequations} \label{OPF}
\begin{align}
\min_{\theta,p^G} \; &c(p^G)=\sum_i(a_{1,i}p_i^G+a_{2,i}(p_i^G)^2) \label{OPFObjective}  \\ 
\mathrm{s.t.}\;&\mathbf{\emph{p}}_i^G-v_i\sum_jv_jB_{ij}^{\text{Pre}}\sin\theta_{ij} = 0  \label{PF} \\
   &\underline{v} \le v_i \le \overline{v} \label{bounds} \\
    &\underline{\theta}_{ij} \le  \theta_{ij}\le  \overline{\theta}_{ij} \label{bounds} 
\end{align}
\end{subequations}
where $p_i^G$ is the power output of generator $i$; $\theta_i$ and $v_i$ are the steady-state phase angle and voltage magnitude of node $i$; and $\theta_{ij}=\theta_i-\theta_j$; $B^\text{Pre}_{ij}$ is an entry of the susceptance matrix $B^\text{Pre}$ of the pre-fault network; (\ref{OPFObjective}) is the total cost function of power generation; constraint (\ref{PF}) represents the steady-state power flow. In order to reveal the conflict between cost-efficiency and stability, we do not consider active and reactive generation bounds for generators in the above power flow model such that the reactive power flow constraints can be eliminated. 

In this three-bus system, bus 1 is the reference bus, so that is phase angle is assumed to be zero. The loads are $p_2^L=1.2$ p.u. and $p_2^L=0.0378$ p.u. The pre-fault susceptance matrix is
$$B^\text{Pre}=\left[
        \begin{array}{ccccc}
          -1.835 \quad & 0.739 \quad & 1.096\\
         0.739 \quad & -1.584 \quad & 0.845 \\
          1.096 \quad & 0.845 \quad & -1.941
        \end{array}
      \right].$$
The optimal solution of (\ref{OPF}) is tabulated in Table \ref{Table1}. 
\begin{table}[h]
\centering
\caption{Optimal Solution of OPF (\ref{OPF})}
\label{Table1}
\begin{tabular}{ccccc}
\hline\hline
\textbf{Node} & \textbf{\begin{tabular}[c]{@{}c@{}}Generation\\ (MW)\end{tabular}} & \textbf{\begin{tabular}[c]{@{}c@{}}Angle\\ (degree)\end{tabular}} & \textbf{\begin{tabular}[c]{@{}c@{}}Objective\\ value (\$)\end{tabular}} & \textbf{Stability}        \\
\hline
1             & 66.91                                                               & 0                                                               & \multirow{3}{*}{2500.4}                                             & \multirow{3}{*}{Unstable} \\
2             & 16.73                                                               & -46.63                                                          &                                                                    &                           \\
3             & 40.14                                   & -8.78                                                          &                                                                    & \\
\hline\hline
\end{tabular}
\end{table}

We adopt the following mathematical model to describe the stability of a power system
\begin{subequations} \label{PSDynamic}
\begin{align}
  &m_i \ddot{\delta_i}+ d_i \dot{\delta_i}  =\mathbf{\emph{p}}_i^G-v_i\sum_jv_jB_{ij}(t)\sin\delta_{ij} \label{Swing} \\
 &\delta(t)\rightarrow \theta \; \text{as}\; t  \rightarrow +\infty, \label{StabilityCertificate}
\end{align}
\end{subequations}
where $\delta_i$ is the dynamic phase angle of bus $i$, where $\delta_i|_{t=t_0}=\theta_i$, and $\delta_{ij}=\delta_i-\delta_j$. The differential equation (\ref{Swing}) of rotor angle dynamics is referred to as "swing" equation in literature\cite{Kundur} while (\ref{StabilityCertificate}) is the criterion for stability. In the above swing equation, an important assumption is that the voltage magnitude $v_i$ for all buses are constant (i.e. $\dot{v}_i=0$ $\forall i$) in the period [$t_0^+,\, t_c^-$]. A short-term outage of Line 12 is considered as the disturbance in this example where the fault results in line tripping. Then, the fault self-clears and Line 12 is re-closed after a short period (e.g. $t_c-t_0=$ 0.1 second), such that :
\begin{align}
    B_{ij}(t)=\begin{cases}
    B^\text{Pre}_{ij} \quad (t \in [t_0^-]\cup [t_c^+,\, +\infty]) \\
    B_{ij}^\text{On} \quad (t\in [t_0^+,\, t_c^-])
    \end{cases} \nonumber \\
    B^\text{On}=\left[
        \begin{array}{ccccc}
          -1.096 \quad & 0 \quad & 1.096\\
         0 \quad & -0.845 \quad & 0.845 \\
          1.096 \quad & 0.845 \quad & -1.941
        \end{array}
      \right]. \nonumber
\end{align}

\begin{figure}[h]
\centering
\includegraphics[width=0.45\textwidth]{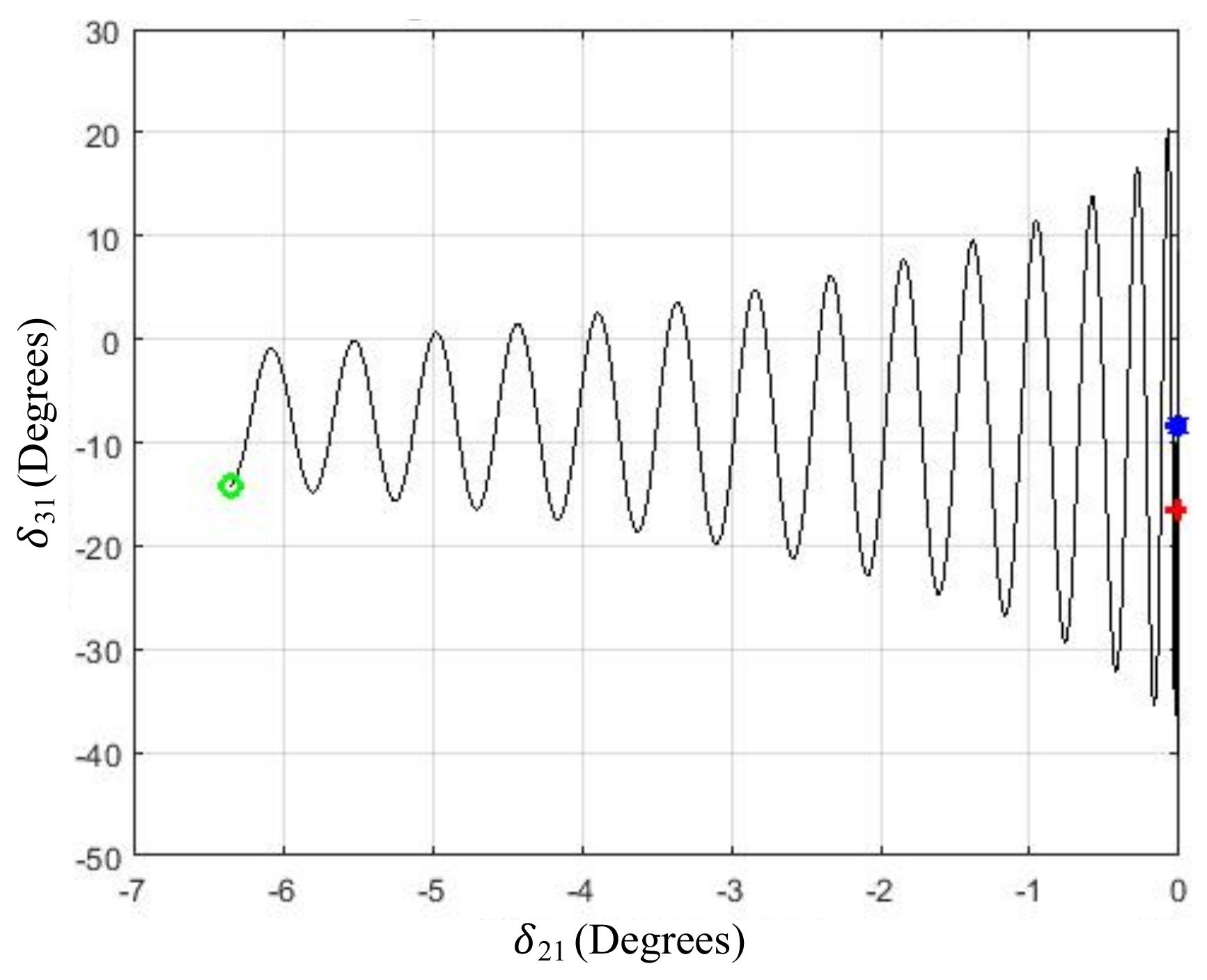}
\caption{System trajectory with the optimal solution of (\ref{OPF}) as the initial condition. The blue asterisk denotes the equilibrium point, the red plus sign represents the fault-cleared state, while the green circle denotes the state of the end of simulation.}
\label{fig:Trajectory1}
\end{figure}

Under this disturbance, the system loses stability when the optimal solution of (\ref{OPF}) tabulated in Table \ref{Table1} is used as the initial condition for dynamic model (\ref{Swing}). The system trajectory is plotted in Figure \ref{fig:Trajectory1}. The dynamics of power system is not considered in the optimization model (\ref{OPF}). As a result, if the optimal solutions of (\ref{OPF}) are used as initial conditions, the stability can not be guaranteed whenever a disturbance occurs. The existing methods \cite{Gan,Scala,Mak,Hijazi} directly impose the stability criterion (\ref{PSDynamic}) into (\ref{OPF}), which results in a computationally intractable large-scale optimization problem. Next subsection applies the proposed approach to develop a computationally-scalable framework for transient stability-constrained optimal power flow (TSCOPF). 

\subsection{Transient Stability-Constrained Optimal Power Flow}

Note that the fault duration ($t_c-t_0$) is relatively small compared with the time constant of dynamic response of system (\ref{Swing}). With the correspondence given in Table \ref{Table2}, it is straight forward to reformulate the TSCOPF problem (\ref{OPF}) in the general expression (\ref{o-SCO}).

\begin{table}[h]
\centering
\caption{Link Between (\ref{OPF}) and (\ref{o-SCO})}
\label{Table2}
\begin{tabular}{cc}
\hline\hline
\begin{tabular}[c]{@{}l@{}}Variables or constraints\\  in Problem (\ref{OPF})\end{tabular} & \begin{tabular}[c]{@{}l@{}}Variables or constraints \\ in Problem (\ref{o-SCO}) \end{tabular} \\
\hline \\                                                                                      $p^G$ & $w$     \\                                                        $\theta$  &  $x^o$ \\
                                                                                       $\delta$     &     $x$         \\
                                                                                        $B(t)$    & $u(t)$       \\
                                                                                         (\ref{PF})   &     (\ref{NonLinearSyst})         \\
                                                                                           (\ref{Swing})     &      (\ref{DynamicSyst})                                                                                      \\
                                               (\ref{bounds})                                             &      (\ref{PhysicalConstraint})                                                                                      \\
                                              (20e)                                              &                                                (\ref{Converge})                                            \\
 \hline\hline                                                                                        
\end{tabular}
\end{table}

It has been shown in \cite{LongA,LongSynchronization} that the dynamic system (\ref{Swing}) is a Lur'e system by considering
\begin{equation}
\phi_i (C(\delta-\theta))=v_i\sum_jv_jB_{ij}(t)(\sin\theta_{ij}-\sin\delta_{ij}) \nonumber
\end{equation}
as the nonlinear term. For the detailed process of reformulating (\ref{Swing}) as a Lur'e-form system, please refer to \cite{LongA,LongSynchronization}. Generally, $-\pi /2 \le \theta_{ij} \le \pi /2\;(\forall ij)$. As a result, we choose $\mathcal{P}=\{\delta \mid -\pi /2 \le \delta_{ij} \le \pi /2,\;\forall ij\}$, and consequently
\begin{equation}
[\gamma,\;\beta]=[\frac{\underline{V}^2\underline{B}(1-\sin|\underline{\theta}|)}{ \pi /2-|\underline{\theta}|},\;\overline{V}^2\overline{B}], \nonumber
\end{equation}
where $\underline{B}$ and $\overline{B}$ are the minimum and maximum elements of $B^\text{Pre}$ respectively; and $\underline{\theta}$ is the minimum angle difference over all the lines. With the chosen sector, LMI (\ref{LMI}) can be solved to obtain a positive definite matrix $P$. Then, we can obtain the following scalable NLP problem by applying the stability-constrained optimization framework (SCO1) to problem (\ref{OPF}):
\begin{align}  \tag{TSCOPF}
\begin{split} \label{TSCOPF}
\min_{\theta,p^G,V^{\min}}  &c=\sum_i(a_{1,i}p_i^G+a_{2,i}(p_i^G)^2)- \epsilon V^{\min}  \\
\text{s.t.} \quad &\text{(\ref{PF}),\;(\ref{bounds})},\;\text{and}  \\
&\delta_i^c=\theta_i+\frac{(t_c-t_0)^2}{2!}K_i \\
&V(\delta_i^c) \le V^{\min}\\
&V^{\min} \le \frac{(\theta_{ij} - \pi /2)^2}{C_i^\mathrm{T}
P^{-1} C_i} \\
&V^{\min} \le \frac{(\theta_{ij} + \pi /2)^2}{C_i^\mathrm{T}
P^{-1} C_i},
\end{split}
\end{align} 
where $K_i=V_i\sum_jV_j(B^\text{On}_{ij}-B^\text{Pre}_{ij})\sin\theta_{ij}$; $C_i$ is a column vector with its $i$th and $j$th elements equal to 1 and -1 respectively, while all the other elements are 0; the expression of $\delta_i^c$ stems from (\ref{trajectory}) with $N=3$. Note that the fault clearing time is generally 0.05-0.1 second. Thus, the terms with $(t_c-t_0)^n$ is negligible if $n>3$. 

\begin{figure}[h]
\centering
\includegraphics[width=0.45\textwidth]{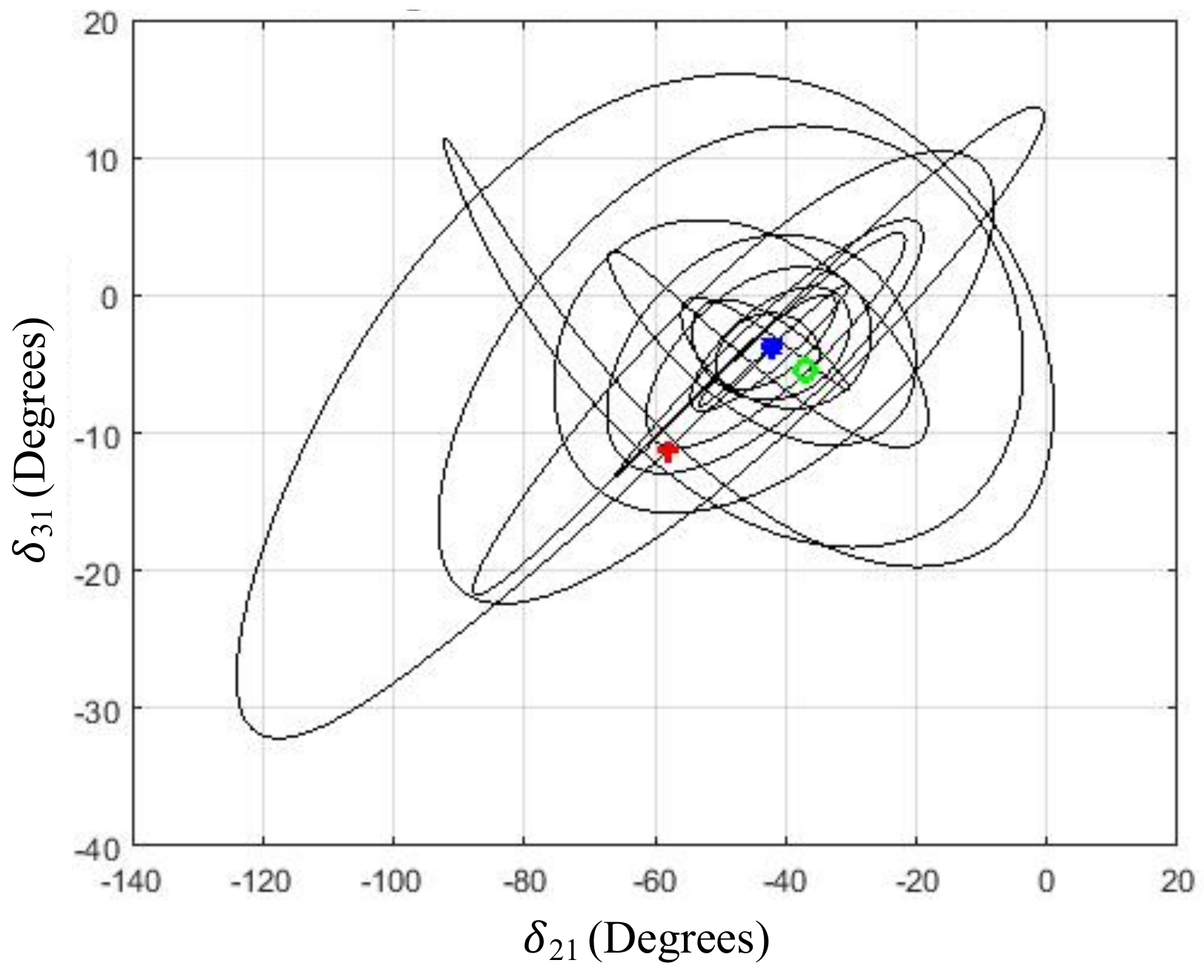}
\caption{System trajectory with the optimal solution of (\ref{TSCOPF}) as the initial condition. The blue asterisk denotes the equilibrium point, the red plus sign represents the fault-cleared state, while the green circle denotes the state of the end of simulation.}
\label{fig:Trajectory2}
\end{figure}

The optimal solution of (\ref{TSCOPF}) is tabulated in Table \ref{Table3}, where the optimal cost is a little bit higher than the optimal cost of (\ref{OPF}). With this optimal solution as the initial condition, the system is stable after the disturbance is cleared. The system trajectory of this case is plotted in Figure \ref{fig:Trajectory2}.
\begin{table}[h]
\centering
\caption{Optimal Solution of OPF (\ref{OPF})}
\label{Table3}
\begin{tabular}{ccccc}
\hline\hline
\textbf{Node} & \textbf{\begin{tabular}[c]{@{}c@{}}Generation\\ (MW)\end{tabular}} & \textbf{\begin{tabular}[c]{@{}c@{}}Angle\\ (degree)\end{tabular}} & \textbf{\begin{tabular}[c]{@{}c@{}}Objective\\ value (\$)\end{tabular}} & \textbf{Stability}        \\
\hline
1             & 55.23                                                               & 0                                                               & \multirow{3}{*}{2501.3}                                             & \multirow{3}{*}{Stable} \\
2             & 20.16                                                               & -42.35                                                          &                                                                    &                           \\
3             & 48.39                                                               & -3.73                                                          &                                                                    & \\
\hline\hline
\end{tabular}
\end{table}

\textbf{Remark 2}. To guarantee a solution that satisfies the stability criterion (\ref{StabilityCertificate}), the existing approach needs to discretize the dynamic constraint (\ref{Swing}) into hundreds or even thousands algebraic equations, which will make the problem computationally intractable. In contrast, the proposed framework (\ref{TSCOPF}) requires only a very limited number of algebraic constraints to replace the differential constraints. That makes (\ref{TSCOPF}) highly-scalable for large-scale power grids which have hundreds to thousands buses. Note that not all the dynamic models of power systems can be formulated in the Lur'e form. The quasi-polynomial system (\ref{Pseudo}) discussed in Subsection IV-B offers an opportunity to apply the proposed approach to solve the stability-constrained optimization for the non-Lur'e cases in power systems.

\section{Conclusion}
This paper proposes a novel method to solve a stability-constrained optimization problem, where the constraints include differentiate-algebraic equations (DAEs). The proposed approach converts the DAE-based stability certificate into a small set of algebraic constraints, leading to a convex optimization framework, based on the LaSalle's invariance principle. Consequently, it is potentially computationally-tractable even for large-scale nonlinear systems offering a viable rigorous alternative to the existing discretize-then-optimize approaches.

Theoretically, the proposed solution method can be applied to an broad range of nonlinear systems, specifically the ones for which there exists a convex common local Lyapunov function (CLLF). This paper discusses the methods of constructing quadratic and SOS-convex CLLFs for Lur'e-type and quasi-polynomial systems respectively. The applications of the proposed approaches in these two types of nonlinear systems are introduced and illustrated via a numerical experiment on a three-generator power system. The simulation results demonstrates the effectiveness of the proposed approach.

\appendices



\section{Formulation of the Modified Pendulum System} \label{app:pendulum}
The modified pendulum system used in the Figure \ref{fig:pulse} is given as
\begin{align}
& \dot{x}_1 = x_2 - u(t) \nonumber
\\
& \dot{x}_2 = -10 \sin (x_1)-x_2 +u(t). \nonumber
\end{align}

\section{Proof of Lemma 1}\label{app:lemma1}

The disjunctive program (\ref{Vmindefinition}) implies that $V^{\min}$ is the minimum value of $V$ alone the facets of polytope $\mathcal{P}$. Suppose $x^*$ is a point on one of the facets of the  $\mathcal{P}$, we have $V(x^*) \geq V^{\min}$. Hence, the system is not able to evolve from the fault clearing point $x^c$ to $x^*$ since the Lyapunov function value $V(x)$ can only decrease along the system trajectory since $\dot{V} \le 0$ within the set $\mathcal{P}$ and $V(x^c) \le V^{\min} \le V(x^*)$. As a result, the system trajectory will stay within the region $\Omega$ or even converge to the origin as $t \rightarrow +\infty$.

\section{Proof of Theorem 1} \label{app:theorem1}
For each $i \in \mathcal{N}$, constraints (5b) and (5d) contain $n$ equations and $|\mathcal{N}|$ inequalities respectively, while (5c) contains only one equation. Constraints (5b)-(5d) specify a unique solution $\hat{x}^i$ to $x^i$ as the number of variables is ($n+1$). As a result, we can let $\hat{X}=[\hat{x}^1,\ldots,\, \hat{x}^k,\ldots,\, \hat{x}^{|\mathcal{N}|}]$ and $\hat{V}^{\min}=V(\hat{x}^k)$ denote the optimal solution of (\ref{newVmin}).

The DP problem (\ref{Vmindefinition}) is equivalent to the following bilevel NLP problem
\begin{align}
V^{\min}=\min_{i\in \mathcal{N}}\; &V_i^{\min}  \nonumber\\
\begin{split} \label{subproblem}
 V_i^{\min}=\min_{x}\; &V(x) \\
\mathrm{s.t.}\; &C_i^{\mathrm{T}}x -d_i= 0 \\
&Cx-d \le 0 .
\end{split}
\end{align}
It is sufficient to show that $\hat{x}^i$ is optimal to the $i$th subproblem in (\ref{subproblem}) since $\hat{x}^i$ satisfies the first-order \textit{K-K-T} conditions (5b)-(5d) and $V$ is convex. According to the objective function (5a), it suffices to show that $\hat{x}^k$ is the optimal solution to the main problem of (\ref{subproblem}). Hence, the $V^{\min}$ obtained by solving problem (\ref{newVmin}) is the same as that obtained by solving (\ref{Vmindefinition}).

\section{Proof of Theorem 2} \label{app:theorem2}

Assume that $\hat{s}=[\hat{x}^o$, $\hat{x}^1,\ldots, \hat{x}^{|\mathcal{N}|}$, $\hat{w}$, $\hat{V}^{\min}]$ is an optimal solution of (s-SCO):\\
\textbf{i.} $\hat{s}$ is feasible to (b-SCO). 

First, we will show that any optimal solution of (s-SCO) is optimal to problem (\ref{newVmin}). Suppose that $\hat{V}^{\min}$ does not make equal sign hold in any of (\ref{Vminbound}), which means $\hat{s}$ is not feasible to (\ref{newVmin}). It suffices to show there exists another feasible solution of (s-SCO), $\bar{s}=[\hat{x}^o$, $\hat{x}^1,\ldots, \hat{x}^{|\mathcal{N}|}$, $\hat{w}$, $\hat{V}^{\min}+\Delta V]$, where $\Delta V$ is an arbitrarily small positive value. We have  
\begin{equation}
c^\prime(\bar{s})-c^\prime(\hat{s})=-\epsilon \Delta V \le 0, \nonumber
\end{equation}
which contradicts the optimality of $\hat{s}$. In other words, when substitute the optimal solution $\hat{s}$ of (s-SCO) into inequalities (\ref{Vminbound}), equal sign holds in at least one of them. It implies that $\hat{s}$ is optimal to problem (\ref{newVmin}) and, consequently, feasible to problem (b-SCO), since optimization problem (\ref{newVmin}) is a constraint of (b-SCO).\\
\textbf{ii.} $\hat{s}$ is also optimal to (b-SCO). 

Suppose $\hat{s}$ is not a locally optimal solution of (b-SCO), then there exist a feasible solution of (b-SCO), $\tilde{s}=[\tilde{x}^o$, $\tilde{x}^1,\ldots, \tilde{x}^{|\mathcal{N}|}$, $\tilde{w}$, $\tilde{V}^{\min}]$, that is in the vicinity of $\hat{s}$ satisfying
\begin{equation}
c(\tilde{x}^o, \tilde{w})- c(\hat{x}^o,\hat{w}) \le 0. \nonumber
\end{equation}
It is straightforward to show that all solutions which are feasible to subproblem (\ref{newVmin}) will also satisfy constraint (\ref{Vminbound}). Thus, $\tilde{s}$ is also feasible to (s-SCO) and satisfies
\begin{equation} \label{Contradition}
c^\prime(\tilde{s})-c^\prime(\hat{s})=c(\tilde{x}^o, \tilde{w})- c(\hat{x}^o,\hat{w}) + \epsilon (\hat{V}^{\min}-\tilde{V}^{\min}) \le 0.
\end{equation}
Note that $\epsilon$ is an arbitrarily small value. Hence, it is reasonable to assume that the term $\epsilon$($\hat{V}^{\min}$ - $\tilde{V}^{\min}$) is not comparable to ($c$($\tilde{x}^o, \tilde{w}$) - $c$($\hat{x}^o,\hat{w}$)), which means ($c^\prime$($\tilde{s}$) - $c^\prime$($\hat{s}$)) has the same sign as ($c$($\tilde{x}^o, \tilde{w}$) - $c$($\hat{x}^o,\hat{w}$)). Condition (\ref{Contradition}) contradicts the assumption that $\hat{s}$ is not locally optimal to (b-SCO). Namely, $\hat{s}$ is an local minimum of (b-SCO) if it is an local minimum of (s-SCO).

So far, the Theorem 2 has been proved by contradiction.

\section{Proof of Theorem 3} \label{app:theorem3}

For the sake of convenience, we replace the terms $C_i^\mathrm{T}
x^o$ and $C_i^\mathrm{T}
P^{-1} C_i$ with $X$ and $1/\mu$ respectively. The notation $CONV(A)$ means the convex hull of set $A$.\\
\textbf{i.} $CONV(\psi) \subseteq \Psi$.

For any $X \in [\underline{X},\;(\overline{d}+\underline{d})/2]$, we have 
\begin{align}
\psi &=\{(X,V^{\min})|0 \le V^{\min} \le \mu(X-\underline{d})^2\} \nonumber \\
\Psi &=\{(X,V^{\min})|0 \le V^{\min} \le \mu(a_{i}^\prime X-b_i^\prime)\}. \nonumber
\end{align}
It suffices to know that $\psi \subseteq \Psi$ since $(X-\underline{d})^2 \le (a_{i}^\prime X-b_i^\prime)$ in the interval [$\underline{X},\;(\overline{d}+\underline{d})/2$]. Similarly, for any $X \in [(\overline{d}+\underline{d})/2,\; \overline{X}]$, we have the same conclusion that $\psi \subseteq \Psi$, which means $\Psi$ is a convex relaxation of $\psi$. Since convex hull is defined as the intersection of all convex relaxations of a non-convex set \cite{LiConvex2,LiThe}, we have $CONV(\psi) \subseteq \Psi$.\\
\textbf{ii.} $CONV(\psi) \supseteq \Psi$.

If a linear inequality is valid for a given set $\Omega_A$, it will also be valid for any subset of $\Omega_A$. Note that “a linear inequality is valid for a set” means the inequality is satisfied by all its feasible solutions \cite{Validcut}. On the other way round, according to the properties of supporting hyperplanes \cite{BoydConvex}, $\Omega_B$ is said to be a subset of $\Omega_A$ if $\Omega_A$ is convex and any valid linear inequality of $\Omega_A$ is also valid for $\Omega_B$ \cite{LiConvex2}. Let $s=[X$, $V^{\min}]$ and suppose that $\alpha s \geq \beta$ is any valid linear cut for $CONV(\psi)$, this cut should be also valid for all the points in $\psi$. To prove that $CONV(\psi) \supseteq \Psi$ (i.e. $\Psi$ is a subset of $CONV(\psi)$), we try to show that $\alpha s \geq \beta$ is valid for all the edges of $\Psi$.

The convex set $\Psi$ has five edges of which the formulations are given as
   \[   
    \Psi_1 = \left\{(X,V^{\min}) \left| \begin{array}{lr}
V^{\min}= \lambda ((\overline{X} - 2\Delta l)X+\Delta l^2) \\ 
W^{\min} \le \lambda ((\underline{X} + 2\Delta l)X+\Delta l^2) \\
\underline{X} \le X \le \overline{X}
\end{array}\right. \right\}, 
  \]
     \[   
    \Psi_2 = \left\{(X,W^{\min}) \left| \begin{array}{lr}
W^{\min} \le \lambda ((\overline{X} - 2\Delta l)X+\Delta l^2) \\ 
W^{\min} = \lambda ((\underline{X} + 2\Delta l)X+\Delta l^2) \\
\underline{X} \le X \le \overline{X}
\end{array}\right. \right\}, 
  \]
     \[   
    \Psi_3 = \left\{(X,W^{\min}) \left| \begin{array}{lr}
W^{\min} \le \lambda ((\overline{X} - 2\Delta l)X+\Delta l^2) \\ 
W^{\min} \le \lambda ((\underline{X} + 2\Delta l)X+\Delta l^2) \\
 X = \underline{X}
\end{array}\right. \right\}, 
  \]
     \[   
    \Psi_4 = \left\{(X,W^{\min}) \left| \begin{array}{lr}
W^{\min} \le \lambda ((\overline{X} - 2\Delta l)X+\Delta l^2) \\ 
W^{\min} \le \lambda ((\underline{X} + 2\Delta l)X+\Delta l^2) \\
X = \overline{X}
\end{array}\right. \right\}, 
  \]
     \[   
    \Psi_5 = \left\{(X,W^{\min}) \left| \begin{array}{lr}
W^{\min} = 0 \\
\underline{X} \le X \le \overline{X}
\end{array}\right. \right\}. 
  \]
  
As an example, we show that the cut $\alpha s \geq \beta$ is valid for edge $\Psi_1$ in this paragraph. It is easy to verify that the two points $s_1=(0,\lambda \Delta l^2)$ and $s_2=(\overline{X},\lambda(\overline{X}-\Delta l)^2)$ are located in both $\psi$ and $\Psi_1$. That means the cut $\alpha s \geq \beta$ is valid for these two points and we have $\alpha s_1 \geq \beta$ and $\alpha s_2 \geq \beta$. Let $\hat{s}=(\hat{X},\hat{W}^{min})$ denote any given point in set $\Psi_1$. For any given $\hat{s}$, there exists a value $c$ $(0 \le c \le 1)$ satisfying $\hat{s}=cs_1+(1-c)s_2)$. It suffices to verify this statement by substituting $\hat{s}=((1-c)\overline{X},c\lambda \Delta l^2 +(1-c)\lambda(\overline{X}-\Delta l)^2)$ into the first equation in $\Psi_1$. As a result, we have
\begin{equation}
\alpha \hat{s}= c\alpha s_1+(1-c) \alpha s_2 \geq c\beta +(1-c)\beta = \beta, \nonumber
\end{equation}
which means the linear cut is also valid for any given point in $\Psi_1$ and, consequently, valid for $\Psi_1$.

Using the same method, the readers are able to prove that the linear cut $\alpha s \geq \beta$ is valid for all the other four edges of $\Psi$ and, consequently, valid for the whole convex set $\Psi$. Hence, $CONV(\psi) \supseteq \Psi$.\\
\textbf{iii.} $\Psi = CONV(\psi)$.

$CONV(\psi) \subseteq \Psi$ and $CONV(\psi) \supseteq \Psi$ together imply $\Psi = CONV(\psi)$.

\section{Application in Polynomial Systems} \label{app:polynomial}

If the dynamic system (\ref{nominal}) is a polynomial system, it can be rewritten as
\begin{equation}  \label{Polynomial}
\dot{x}=f(x-x^o)=[A+G(x-x^o)](x-x^o)
\end{equation}
where $G(\cdot):\mathbb{R}^n \rightarrow \mathbb{R}^n$ $(i=1,\ldots,n)$, $G(0)=0$,  and matrix $A$ is hurwitz. Let
\begin{equation}
\mathcal{P}=\{x \mid \|x-x^o\|\le r\},\nonumber
\end{equation}
where $r$ is some given value.

\textbf{Proposition 4} \textit{The quadratic function} (\ref{QuaLyapunov}) \textit{is a valid CLLF of the nonlinear system} (\ref{Polynomial}) \textit{if the positive definite matrix $P$ satisfies, $\forall (x^o,w) \in \mathcal{F}$} \\
\textit{a) $PA+A^TP=-Q$ is negative definite,}\\ \textit{b)} $2 \gamma \|P\|\le -\lambda_{\min}(Q),$\\
\textit{where $\lambda_{\min}(\cdot)$ denotes the minimum eigenvalue of a matrix, and}
\begin{align}
\begin{split}
\gamma=\max_{x,w,x^o}\ &\|G(x-x^o)\|  \\
\mathrm{s.t.}\ &x \in \mathcal{P}\  \text{and}\  (w,\ x^o) \in \mathcal{F}.\nonumber
\end{split}
\end{align}

\textit{Proof}: Let $y=(x-x^o)$For any given ($x^o,w$) $\in \mathcal{F}$, the derivative of quadratic CLLF $V$ along the trajectories of the nonlinear system in (\ref{Polynomial}) is given by
\begin{align} \label{Vdot}
\dot{V}(y)&=\dot{y}^TPy+y^TP\dot{y} \nonumber \\
&=y^T(PA+ A^TP)y+2y^TPG(y)y \nonumber \\
&\le -\lambda_{\min}(Q)\|y\|^2+2\|P\|\|G(y)\|\|y\|^2. \nonumber
\end{align}
According to condition \textit{b)} in the proposition, it suffices to show that $\dot{V}<0$ for all $(w,\ x^o) \in \mathcal{F}$. $\square$

Applying the single-level stability-constrained optimization framework (s-SCO) to the polynomial system cases (\ref{Polynomial}), we have
\begin{align}  
\begin{split}
\min_{z} \quad &\text{(\ref{Objective})} \\
\text{s.t.} \quad &\text{(\ref{NonLinearSyst}),\,(\ref{PhysicalConstraint})},\,(\ref{trajectory}),\,\text{and}\\
&V(x^c) \le \lambda_{\min}(P)r^2,
\end{split} \tag{SCO3}
\end{align}
since
\begin{equation}
\min_{ \|x-x^o\|= r}V(x)=\lambda_{\min}(P)r^2. \nonumber
\end{equation}

\end{document}